# A Modified Burg Algorithm Equivalent In Results to Levinson Algorithm

Rami Kanhouche

*Abstract*— We present a new modified Burg-Like algorithm for spectral estimation and adaptive signal processing that yield the same prediction coefficients given by the Levinson algorithm for the solution of the normal equations. An equivalency proof is given for both the 1D signal and 2D signal cases. Numerical simulations illustrate the improved accuracy and stability in spectral power amplitude and localization; especially in the cases of low signal to noise ratio, and (or) augmenting the used prediction coefficients number for a relatively short data records. Also our simulations illustrate that for relatively short data records the unmodified version of Burg Algorithm fail to minimize the mean square residual error beyond certain Order, while the new algorithm continue the minimization with Order elevation.

## I. INTRODUCTION

So many are the fields of application of Linear Autoregressive Modelling in our days [1]-[3]. At the same time, a lot of work has been done in studying this model and its applications in the one-dimensional signal case, and -lately- for the two-dimensional signal. Almost One common factor to the previous work done depend on the solution of the normal equations –known also as the Yule-Walker equations- using different techniques to find an estimation of the prediction coefficients, which are scalars in 1D case, and matrices in the Multichannel (multivariate) case. Starting from a given 1D signal, two main methods of estimating the prediction coefficients are the Levinson Algorithm [7] and the Burg Algorithm [6]. The counterparts of these two algorithms in the 2D case are, respectively, a modified version of the Whittle-Wiggins-Robinson Algorithm (WWRA)[8][9], and a group of slightly in-between different algorithms given in [10], [11], [15], which will be referred to in this paper–excluding [10], as 2D burg Algorithms. Historically the 2D approach depended on the Multichannel Levinson (WWRA) and Multichannel Burg [12] approach, to obtain the Coefficients matrices, with an extra step to recover the 2D quarter-plan causal Filter scalar coefficients, from the calculated coefficients Matrices, the connection between the Coefficients Matrices and the prediction 2D Filter was first established by [14], this was further invested by [11], to apply the Multichannel burg algorithm [12], directly on the two dimensional data, which gave the first 2D-Burg algorithm, the same algorithm was represented in [15], with an additional simplifications that take advantage of the Toeplitz-Block Toeplitz of the 2D Correlation Matrix, while it is just Block Toeplitz in the Multichannel case. On the other side in [10], another Burg-like algorithm was presented, that calculate the 2D Prediction Filter Directly without passing by the calculus of the corresponding Multichannel Coefficients Matrices.

In this paper we present a new 1D and 2D Burg Algorithm, that yield the same prediction coefficients as the Levinson Algorithm in 1D case, and the WWRA in 2D case. Also we enhance WWRA computational efficiency in 2D case, due to proposition 3. In practice, the motivation of using 1D Burg is the computational efficiency over the Levinson Algorithm, since the later, demand the calculus of the correlation matrix in advance, which is very costly in calculus time. For the same reason also the 2D Burg Algorithms are preferred over the WWRA algorithm. Still, for the reasons that we will explain in this paper, the 1D and 2D Burg Algorithms suffer from major numerical deficiency, that was the motive behind our development of the new Burg Algorithm. While the computational efficiency of The New Modified Burg Algorithm need to be further investigated they yield a stable Levinson and WWRA solution. Saying that, we also believe, that the new-presented Burg Algorithm, in 1D and 2D cases, represent the *natural burg algorithms*, since they are now, consistent, in theory, and in application, to the original Levinson approach.

In [4]-[5], fast, and practical formulas were presented, that calculate the inverse of the Correlation Matrix in function of the prediction coefficients.

The paper is organized as follows. In section II, we establish the proof of asymptotic equivalency, between, Levinson Algorithm results, and Burg Algorithm results, in 1D case. In section III, we also establish the proof, of asymptotic equivalency, between WWRA Algorithm results, and the 2D Burg Algorithms results. In section IV, we present the practical modifications, needed to the current Burg Algorithms, in both, 1D and 2D cases, to realize the *in-results* proven equivalency, to Levinson algorithm, and WWRA algorithm, respectively. Finally, Section V contains our

Manuscript received September 23, 2003.
The Author is with Lab. CMLA, Ecole Normale Supérieure, 61, avenue du Président Wilson, 94235 CACHAN Cedex, France.
phone: +33-1-64971240;fax: +33-1-546546546;
e-mail: rami.kanhouche@cmla.ens-cachan.fr.

numerical examples, which illustrate, the numerical stability of the new Burg algorithms, relative to the original Burg Algorithms.

## II. 1D CASE

### A. The model

The autoregressive Model

$$x(k) + \sum_{l=1}^{n} a_l^n x(k-l) = w(k) \qquad (1)$$

where $w$ represent the prediction error of the signal $x$, yield the Normal Equations

$$R_n \mathbf{a}_n = -\mathbf{r}_n \qquad (2)$$

where $R_n$ is a Toeplitz Hermetian matrix with size (n,n), defined as

$$R_n \equiv [r_{ij}]_{i,j=0,1,..n-1}, \quad r_{ij} = r_{ji}^* = r_{i-j} \qquad (3)$$

$$r_t \equiv E(x(k+t)x(k)^*), t = 0, \pm 1, ..., \pm n \qquad (4)$$

$$\mathbf{a}_n \equiv [a_1^n, a_2^n, ..a_n^n]^T, \mathbf{r}_n \equiv [r_1, r_2, ..r_n]^T \qquad (5)$$

### B. The 1D Levinson Algorithm

The Levinson Algorithm [7] proceed to the solution of (2) in a recursive manner, for $n > 1$

$$a_n^n = -\frac{r_n + \hat{\mathbf{r}}_{n-1}^T \mathbf{a}_{n-1}}{r_0 + \mathbf{r}_{n-1}^T \mathbf{a}_{n-1}^*} \qquad (6)$$

$$\mathbf{a}_n = \begin{bmatrix} \mathbf{a}_{n-1} \\ 0 \end{bmatrix} + a_n \begin{bmatrix} \hat{\mathbf{a}}_{n-1}^* \\ 1 \end{bmatrix} \qquad (7)$$

starting from

$$a_1^1 = -\frac{r_1}{r_0}. \qquad (8)$$

where the notation $\hat{\mathbf{x}}$ define a vector with the same elements of $\mathbf{x}$ taken in reverse order.

### C. 1D Burg Algorithm

The Burg Algorithm [6] define an estimate of the prediction errors

$$e_n^f(k) \equiv \sum_{l=1}^{n} a_l^n x(k-l) + x(k)$$

$$e_n^b(k) \equiv \sum_{l=1}^{n} a_l^{n*} x(k+l-n) + x(k-n) \qquad (9)$$

The prediction coefficient is also calculated in a recurrent manner, according to the relation

$$a_n^n = \frac{-E(e_{n-1}^f(k) e_{n-1}^{b*}(k-1))}{\frac{1}{2} E(|e_{n-1}^f(k)|^2 + |e_{n-1}^b(k-1)|^2)}. \qquad (10)$$

Using the Levinson recurrence relation (7), a recursive relation between the errors signals is obtained and used in the calculus, to update the error signals for the next order prediction coefficient calculus

$$e_{n+1}^f(k) = e_n^f(k) + a_{n+1}^{n+1} e_n^b(k-1) \qquad (11)$$

$$e_{n+1}^b(k) = e_n^b(k-1) + (a_{n+1}^{n+1})^* e_n^f(k). \qquad (12)$$

### D. Equivalency proof

To proceed in our construction, we will need the following proposition, which is easy to verify.

**Proposition.** *The forward and backward prediction errors square means are asymptotically equal,*

$$\lim_{k \to +\infty} \sum_k |e_n^f(k)|^2 = \lim_{k \to +\infty} \sum_k |e_n^b(k-1)|^2. \qquad (13)$$

**Theorem.** *Asymptotically, the prediction coefficients calculated using Levinson Algorithm and Burg Algorithm are Equals.*

Proof. Using the (13) we can know write the Burg relation (10) in the form

$$a_n^n = \frac{-E(e_{n-1}^f(k) e_{n-1}^{b*}(k-1))}{E(|e_{n-1}^b(k)|^2)}. \qquad (14)$$

Replacing (9) into (14), we proceed according to the following

$$a_{n+1}^{n+1} = E \begin{pmatrix} x(k)x(k-n-1) + \\ \sum_{l=1}^{n} x(k)x^*(k-1-n+l)a_l + \\ \sum_{l=1}^{n} a_l^n x(k-l)x^*(k-n-1) + \\ \sum_{l=1}^{n} a_l^n \sum_{l'=1}^{n} x(k-l)x^*(k-1-n+l')a_{l'}^n \end{pmatrix} /$$

$$E \begin{pmatrix} x(k-n)x(k-n) + \sum_{l=1}^{n} x(k-n)x^*(k-n+l)a_l^n \\ + \sum_{l=1}^{n} a_l^{n*} x(k-n+l)x^*(k-n) \\ + \sum_{l=1}^{n} a_l^{n*} \sum_{l'=1}^{n} x(k-n+l)x^*(k-n+l')a_{l'}^n \end{pmatrix} \qquad (15)$$

$$a_{n+1}^{n+1} = -\left( r_{n+1} + \sum_{l=1}^{n} r_{n+1-l} a_l^n \atop + \sum_{l=1}^{n} a_l^n r_{n+1-l} + \sum_{l=1}^{n} a_l^n \sum_{l'=1}^{n} r_{n+1-l-l'} a_{l'}^n \right) / \left( r_0 + \sum_{l=1}^{n} r_{-l} a_l^n \atop + \sum_{l=1}^{n} a_l^{n*} r_l + \sum_{l=1}^{n} a_l^{n*} \sum_{l'=1}^{n} r_{l-l'} a_{l'}^n \right)$$ (16)

Using the equality

$$\sum_{l=1}^{n} a_l r_{d-l} = -r_d \quad d = 1, 2, ..n$$ (17)

which is just another form for the relation (2), the second and forth terms, in both the dominator and denominator of (16) are eliminated

$$a_{n+1}^{n+1} = -\left( r_{n+1} + \sum_{l=1}^{n} r_{n+1-l} a_l \right) / \left( r_0 + \sum_{l=1}^{n} a_l^* r_l \right).$$ (18)

This completes the proof.

## III. 2D CASE

### A. The Model

The Multichannel signal linear Autoregressive model is defined as [16]

$$X(k) = -\sum_{l=1}^{n_1} A_l^{n_1} X(k-l) + W(k)$$ (19)

where $X, W$ are vectors of size $n_2 + 1$, changing according to the time $k$, $A_l^{n_1}$ is the $l$ prediction coefficient Matrix of degree $n_1$ with size $((n_2+1) \times (n_2+1))$.

This Model yield also the multi-channel normal equations, known as multi-channel Yule-walker equations

$$\mathbf{A}_{n_1,n_2} . R_{n_1,n_2} = \mathbf{r}_{n_1,n_2}$$ (20)

$$\mathbf{A}_{n_1,n_2} \equiv \left[ A_1^{n_1}, A_2^{n_1}, ..., A_{n_1}^{n_1} \right]$$ (21)

$$R_{n_1,n_2} \equiv \begin{bmatrix} R_0 & R_1 & \cdots & R_{n_1-1} \\ R_{-1} & R_0 & \ddots & \vdots \\ \vdots & \ddots & \ddots & R_1 \\ R_{1-n_1} & \cdots & R_{-1} & R_0 \end{bmatrix}$$ (22)

$$R_k \equiv E(X(n+k) X^H(n)), k = 0, \pm 1, ..., \pm n_1$$ (23)

$$\mathbf{r}_{n_1,n_2} \equiv [R_1, R_2, \cdots, R_{n_1}]$$ (24)

The Correlation Matrix $R_{n_1,n_2}$ is Hermetian and Block Toeplitz. In the Case of 2D signal the same Multichannel model can be used, by redefining the multi-channel vector, using the 2D signal [14][15]

$$X(k) = \begin{bmatrix} x(k, t_2' + n_2) & x(k, t_2' + n_2 + 1) & \cdots & x(k, t_2'') \\ x(k, t_2' + n_2 - 1) & x(k, t_2' + n_2) & \cdots & x(k, t_2'' - 1) \\ \vdots & \vdots & & \vdots \\ x(k, t_2') & x(k, t_2' + 1) & \cdots & x(k, t_2'' - n_2) \end{bmatrix}$$

$$x(t_1, t_2) \in C, t_1 \in \mathbf{N}, t_2 \in [t_2', t_2''].$$ (25)

In this case the definitions, and relations (19)-(23) still hold for the expression of the Normal equations and the Correlation matrix in the case of the 2D signal, with the Hermetian correlation matrix having the structure Toeplitz-Block-Toeplitz (TBT).

### B. The 2D Levinson Algorithm (WWRA)

The Coefficients Matrix are calculated in a recursive manner, starting from $n > 0$, according to [8][9][16]:

$$A_n^n = -\Delta_n (P_{n-1})^{-1}$$ (26)

$$\Delta_n \equiv R_n + \sum_{l=1}^{n-1} A_l^{n-1} R_{n-l}$$ (27)

$$P_{n-1} \equiv R_0 + \sum_{l=1}^{n-1} J A_l^{n-1*} J R_l$$ (28)

with,

$$\mathbf{A}_{n,n_2} = [\mathbf{A}_{n-1,n_2}, \mathbf{0}] + A_n^n \left[ J A_{n-1}^{n-1*} J, J A_{n-2}^{n-1*} J, \cdots J A_1^{n-1*} J, \mathbf{1} \right]$$ (29)

where $J$ is the exchange Matrix defined as:

$$J \equiv \begin{bmatrix} 0 & \cdots & 0 & 1 \\ \vdots & \ddots & 1 & 0 \\ 0 & \ddots & \ddots & \vdots \\ 1 & 0 & \cdots & 0 \end{bmatrix}$$ (30)

starting From

$$P_0 = R_0$$ (31)

### C. The 2D Burg Algorithm

As in the 1D case, two Prediction Error signals are defined as [11][12][15]

$$\mathbf{e}_n^f(k) \equiv \sum_{l=1}^{n} A_l^n X(k-l) + X(k)$$ (32)

$$\mathbf{e}_n^b(k) \equiv \sum_{l=1}^{n} J(A_l^n)^* J X(k-n+l) + X(k-n)$$ (33)

$n \geq 0$

with recurrent relation

$$\mathbf{e}_{n+1}^f(k) = \mathbf{e}_n^f(k) + A_{n+1}^{n+1} \mathbf{e}_n^b(k-1)$$ (34)

$$\mathbf{e}_{n+1}^b(k) = \mathbf{e}_n^b(k-1) + J\left(A_{n+1}^{n+1}\right)^* J\mathbf{e}_n^f(k) \tag{35}$$

The Prediction Coefficient Matrix of order $n$, is calculated according to [15]

$$A_{n+1}^{n+1} = -\left[P_n^{fb} + J\left(P_n^{fb}\right)^T J\right]\left[P_n^b + J\left(P_n^f\right)^* J\right]^{-1} \tag{36}$$

where

$$P_n^{fb} \equiv E\left(\mathbf{e}_n^f(k)\left(\mathbf{e}_n^b(k-1)\right)^H\right) \tag{37}$$

$$P_n^f \equiv E\left(\mathbf{e}_n^f(k)\left(\mathbf{e}_n^f(k)\right)^H\right) \tag{38}$$

$$P_n^b \equiv E\left(\mathbf{e}_n^b(k)\left(\mathbf{e}_n^b(k)\right)^H\right). \tag{39}$$

The relation (36) result from minimizing the value

$$tr\left(P_n^f + J\left(P_n^b\right)^* J\right). \tag{40}$$

### D. Equivalency proof

**Proposition 2**

$$P_n^b = J\left(P_n^f\right)^* J \tag{41}$$

**Proposition 3**

$$P_n^{fb} = J\left(P_n^{fb}\right)^T J \tag{42}$$

Proof. Using the notation $A_0^n \equiv 1$, and replacing (34), (35) into (37) we get

$$P_n^{fb} =$$

$$E\left(\sum_{l=0}^n \sum_{l'=0}^n A_l^n X(k-l)(X(k-1-n+l'))^H J\left(A_{l'}^n\right)^T J\right) \tag{43}$$

By taking the Expectation inside, and by definition (23) we get

$$P_n^{fb} = \sum_{l=0}^n \sum_{l'=0}^n A_l^n R_{n+1-l-l'} J\left(A_{l'}^n\right)^T J. \tag{44}$$

This can be written in the form

$$P_n^{fb} = R_{n+1} + \sum_{l=1}^n \sum_{l'=1}^n A_l^n R_{n+1-l-l'} J\left(A_{l'}^n\right)^T J$$
$$+ \sum_{l'=1}^n R_{n+1-l'} J\left(A_{l'}^n\right)^T J + \sum_{l=1}^n A_l^n R_{n+1-l} \tag{45}$$

By applying the equality

$$\sum_{l=1}^n A_l^n R_{d-l} = -R_d \quad d = 1,2,..n \tag{46}$$

the second term of (45), can be written as

$$\sum_{l'=1}^n \left(\sum_{l=1}^n A_l^n R_{(n+1-l')-l}\right) J\left(A_{l'}^n\right)^T J = -\sum_{l'=1}^n R_{n+1-l'} J\left(A_{l'}^n\right)^T J \tag{47}$$

which yield

$$P_n^{fb} = R_{n+1} + \sum_{l=1}^n A_l^n R_{n+1-l} \tag{48}$$

Again by rewriting the second term of (45) as

$$\sum_{l=1}^n \sum_{l'=1}^n A_l^n R_{n+1-l-l'} J\left(A_{l'}^n\right)^T J = \sum_{l=1}^n A_l^n \sum_{l'=1}^n \left(J\left(A_{l'}^n\right)^* J R_{n+1-l-l'}^H\right)^H \tag{49}$$

$$= \sum_{l=1}^n A_l^n \left(\sum_{l'=1}^n J\left(A_{l'}^n\right)^* J R_{l'-(n+1-l)}\right)^H$$

and using the equality

$$\sum_{l=1}^n J\left(A_l^n\right)^* J R_{l-d} = -R_{-d}, \quad d = 1,...n \tag{50}$$

we find

$$\sum_{l=1}^n \sum_{l'=1}^n A_l^n R_{n+1-l-l'} J\left(A_{l'}^n\right)^T J = -\sum_{l=1}^n A_l^n \left(R_{-(n+1-l)}\right)^H \tag{51}$$

$$= -\sum_{l=1}^n A_l^n R_{n+1-l}$$

which yield

$$P_n^{fb} = R_{n+1} + \sum_{l'=1}^n R_{n+1-l'} J\left(A_{l'}^n\right)^T J. \tag{52}$$

By considering the fact that for any Toeplitz matrix $T$ the following is always true

$$J(T)^T J = T \tag{53}$$

we apply the transformation $J(\ )^T J$ to the expression (52), we get

$$J\left(P_n^{fb}\right)^T J = R_{n+1} + J\left(\sum_{l'=1}^n R_{n+1-l'} J\left(A_{l'}^n\right)^T J\right)^T J. \tag{54}$$

and by Noticing that for any matrix $M$ we have:

$$J(M)^T J = (JMJ)^T \tag{55}$$

we get

$$J\left(P_n^{fb}\right)^T J = R_{n+1} + \sum_{l'=1}^n \left(J R_{n+1-l'} J\left(A_{l'}^n\right)^T JJ\right)^T$$
$$= R_{n+1} + \sum_{l'=1}^n A_{l'}^n J\left(R_{n+1-l'}\right)^T J \tag{56}$$

By this the proof is complete.

**Theorem 2.** *Asymptotically, the prediction coefficients calculated using WWRA Algorithm and 2D Burg Algorithms are Equals.*

Proof. According to (41) and (42), (36) can be expressed as

$$A_{n+1}^{n+1} = -\left[P_n^{fb}\right]\left[P_n^b\right]^{-1}. \tag{57}$$

Replacing according to (34), (35), we get

$$P_n^{fb} = E\left(\begin{pmatrix} X(k) + \sum_{l=1}^n A_l^n X(k-l) \end{pmatrix} \cdot \begin{pmatrix} (X(k-n-1))^H + \\ \sum_{l=1}^n (X(k-n-1+l))^H \left(J\left(A_l^n\right)^* J\right)^H \end{pmatrix}\right) \tag{58}$$

$$P_n^{fb} = E \begin{pmatrix} X(k)(X(k-n-1))^H + \\ \sum_{l=1}^{n} X(k)(X(k-n-1+l))^H \left(J(A_l^n)^* J\right)^H + \\ \sum_{l=1}^{n} A_l^n X(k-l)(X(k-n-1))^H + \\ \sum_{l=1}^{n} A_l^n \sum_{l'=1}^{n} X(k-l)(X(k-n-1+l'))^H \left(J(A_{l'}^n)^* J\right)^H \end{pmatrix} \quad (59)$$

By taking the expectation over each term, by definition (23) we get

$$P_n^{fb} = \begin{pmatrix} R_{n+1} + \sum_{l=1}^{n} R_{n+1-l} \left(J(A_l^n)^* J\right)^H + \\ \sum_{l=1}^{n} A_l^n R_{n-l+1} + \\ \sum_{l=1}^{n} A_l^n \sum_{l'=1}^{n} R_{n-l-l'+1} \left(J(A_{l'}^n)^* J\right)^H \end{pmatrix} \quad (60)$$

Using the equality

$$\sum_{l=1}^{n} A_l^n R_{d-l} = -R_d \quad d = 1,2,..n \quad (61)$$

which is just another form of (20), we obtain

$$P_n^{fb} = R_{n+1} + \sum_{l=1}^{n} A_l^n R_{n-l+1} . \quad (62)$$

Also by the same way we obtain

$$P_n^b = E \begin{pmatrix} X(k-n)(X(k-n))^H + \\ \sum_{l=1}^{n} X(k-n)(X(k-n+l))^H \left(J(A_l^n)^* J\right)^H + \\ \sum_{l=1}^{n} J(A_l^n)^* J X(k-n+l)(X(k-n))^H + \\ \sum_{l=1}^{n} J(A_l^n)^* J ( \\ \sum_{l'=1}^{n} X(k-n+l)(X(k-n+l'))^H \left(J(A_{l'}^n)^* J\right)^H ) \end{pmatrix} \quad (63)$$

$$P_n^b = \begin{pmatrix} R_0 + \sum_{l=1}^{n} R_{-l} \left(J(A_l^n)^* J\right)^H + \\ \sum_{l=1}^{n} J(A_l^n)^* J R_l + \\ \sum_{l=1}^{n} J(A_l^n)^* J \sum_{l'=1}^{n} R_{l-l'} \left(J(A_{l'}^n)^* J\right)^H \end{pmatrix} \quad (64)$$

Also using the equality

$$\sum_{l=1}^{n} J(A_l^n)^* J R_{l-d} = -R_{-d}, \quad d = 1,...n \quad (65)$$

our Proof is completed.

IV. MODIFIED BURG ALGORITHM

Before going into the New Version of Burg Algorithm it is important to precise an important fact about the calculus of the prediction error auto-correlation. In fact in the case of 2D signal a recurrent relation exist in a way that enhance the calculus cost of both Levinson or Burg Algorithm, according to [15][16]

$$P_n^b = P_{n-1}^b + A_n^n \Delta_n^H \quad (66)$$

By this, and according to the previous sections, Numerically we can consider the equivalent steps in each iteration of both Levinson Algorithm and Burg Algorithm, as the steps concerning the calculus of relations 26, 29, 66, while the difference between the two algorithms reside in the way with which the error cross correlation Matrix is calculated; using 27 in Levinson Algorithm, and 34, 35, 37 in Burg Algorithm, the same type of discussion also apply to the 1D case.

After demonstrating the Asymptotic In Results Equivalency between Levinson Algorithms, and Burg Algorithms in both 1D and 2D case, it is now simple to see the needed modifications to realize this equivalency in the case of limited size signal, in fact the principle modification to the existing Burg algorithms reside in the way relations (11), and (12) in the 1D case, (34), and (35) in the 2D case, are used.

A. *1-D case*

For a finite size signal, $x(k), k \in [0, N-1]$, the current version of burg Algorithm, express the relation (11), and (12) in the form of:

$$e_{n+1}^f(k) = e_n^f(k) + a_{n+1}^{n+1} e_{n+1}^b(k-1) \quad (66)$$

$$e_{n+1}^b(k) = e_n^b(k-1) + (a_{n+1}^{n+1})^* e_{n+1}^f(k)$$

$$k \in [n+1, N-1].$$

As a result of this, the size of the new calculated error signals; diminish with each new order. And as result of this the relation (10) is expressed in the form [6]

$$a_n^n = \frac{-\sum_{k=n}^{N-1} e_{n-1}^f(k) e_{n-1}^{b*}(k-1)}{\frac{1}{2}\sum_{k=n}^{N-1} \left(\left|e_{n-1}^f(k)\right|^2 + \left|e_{n-1}^b(k-1)\right|^2\right)} . \quad (67)$$

The Correction to this is to proceed with the calculus until the last non-zero value, so instead of error signal size diminishing we produce an augmentation in the size of both the forward and backward error signals

$$e_n^b(-1) \equiv 0, e_n^f(N+n) \equiv 0 \quad (68)$$

$$e_{n+1}^f(k) = e_n^f(k) + a_{n+1}^{n+1} e_{n+1}^b(k-1)$$

$$e_{n+1}^b(k) = e_n^b(k-1) + (a_{n+1}^{n+1})^* e_{n+1}^f(k)$$

$$k \in [0, N+n] \quad (69)$$

and the corrected form of the relation (10) becomes

$$a_n^n = \frac{-\sum_{k=1}^{N+n-1} e_{n-1}^f(k) e_{m-1}^{b*}(k-1)}{\sum_{k=0}^{N+n-1} \left(\left|e_{n-1}^f(k)\right|^2\right)}. \tag{70}$$

B. *2D Case*

The same discussion of the 1D case is applied to the 2D case also. In the case of finite size signal

$$x(k,t), k = 0 \ldots N_1 - 1, t = 0 \ldots N_2 - 1 \tag{71}$$

according to the new modified algorithm, the relations (34) (35) are applied according to

$$\mathbf{e}_n^b(-1) \equiv 0, \mathbf{e}_n^f(N+n) \equiv 0 \tag{72}$$

$$\mathbf{e}_{n+1}^f(k) = \mathbf{e}_n^f(k) + A_{n+1}^{n+1} \mathbf{e}_n^b(k-1)$$

$$\mathbf{e}_{n+1}^b(k) = \mathbf{e}_n^b(k-1) + J\left(A_{n+1}^{n+1}\right)^* J \mathbf{e}_n^f(k)$$

$$k \in [0, N+n] \tag{73}$$

and the relation (36) is expressed in the form

$$A_{n+1}^{n+1} = -\left[\sum_{k=1}^{N_1+n} \mathbf{e}_n^f(k)\left(\mathbf{e}_n^b(k-1)\right)^H\right].$$

$$\left[\sum_{k=0}^{N_1+n} \mathbf{e}_n^b(k)\left(\mathbf{e}_n^b(k)\right)^H\right]^{-1} \tag{74}$$

Also to achieve the equivalence in the calculus results with 2D Levinson Algorithm, the 2D Signal need to be distributed into the multiple-vectors X Matrix, in a way that guarantee the Toeplitz block Toeplitz structure of the correlation matrix, so for a limited size 2D signal defined in (71), the definition (25) is applied to the signal according to:

$$X(K) \equiv \begin{bmatrix} x(k,0) & \cdots & x(k,N_2-1) & & \\ & \ddots & & \ddots & \\ & & x(k,0) & \cdots & x(k,N_2-1) \end{bmatrix}$$

So that each row of the signal is repeated and shifted into the next line of the matrix, with zero filled elsewhere.

V. NUMERICAL EXAMPLES

It is important to note, that numerically the modified algorithm is totally equivalent in result to Levinson Algorithm, in both 1D and 2D case. Our objective in this section is to compare Levinson Algorithm application to the Original Burg Algorithm toward spectral estimation, and prediction efficiency. We did apply both of algorithms on a 20 sample of a noised one-sinusoid signals with a relative frequency positioned at $f = 0.25$. The noise signals where generated artificially and added to the spectrum, and the signal samples where obtained using IFFT. Also care has been taken so that all generated signals had precisely a signal to noise ratio equal to 30db.

Figure 1, and 2 represent the spectral estimation versus phase change, resulting from applying Original Burg algorithm, and Levinson Algorithm respectively. The sinusoid phase was changed over 100 equal intervals from zero to 360°. For each phase change the same signal was submitted to each of the two algorithms, while the used prediction order equal to 15. We notice in this case the inaccuracy of Burg Algorithm toward phase change, and toward different noise realisations, on the other hand the Levinson Algorithm manifested much more stable results, toward noise change, and a stable pattern toward phase change. In fact the results displayed in Figure 1 are only readable together because of the use of the log function.
Figure 3, and 4 represent the spectral estimation versus Order change, resulting from applying Original Burg algorithm, and Levinson algorithm respectively. The same signal was submitted to both algorithms with consecutive orders from 1 up until 19, we notice that with order change the frequency identification is stable for both algorithms applied to the same signal, on the other hand Burg algorithm in this case suffer from inconsistency in the frequency power value across different orders.
Finally Figures 5, and 6 displays the Mean square residual error in function of the used Order, resulting from applying Original Burg algorithm, and Levinson algorithm respectively. These results are obtained for the same simulations corresponding to Figure 3, and 4. The error signal was calculated according to relation (1) using the prediction coefficients given by each algorithm, with zero values given for samples outside the provided support. We notice that the original burg algorithm fail to minimize the Mean Square Error for order value above 2, while the Levinson Algorithm continue the minimization with order elevation. This comes to illustrate the Original Burg Algorithm handicap when applied to short data records, or what we can call the *border effect*, this effect is due to the application of the recursive relation (7) without applying the Needed Modifications that we explained, failing to apply the modification (68), (69) will result in a mistaken cross correlation value, and it is simple to see that the effect of this will diminish with augmenting the samples number $N$, relatively to a given Order value $n$.

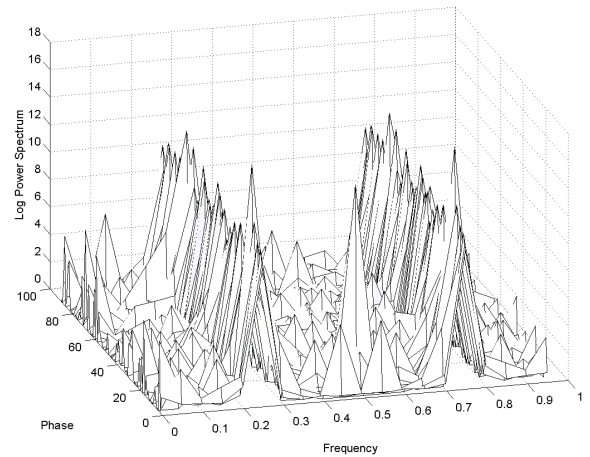

Figure 1

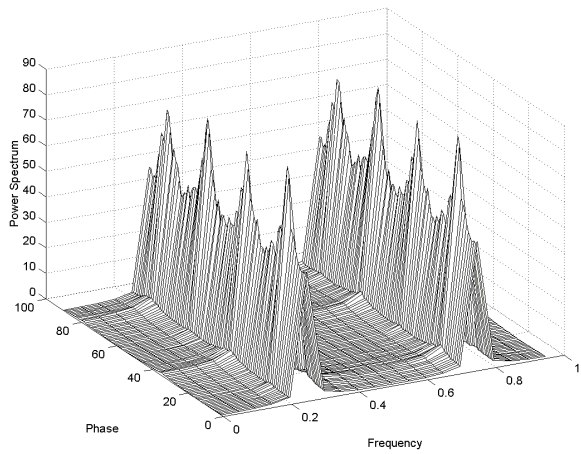

Figure 2

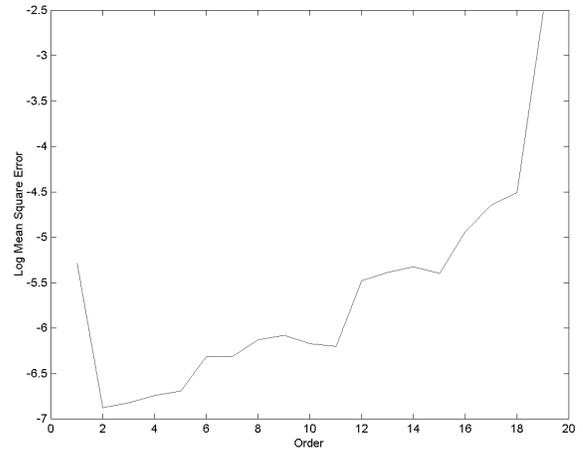

Figure 5

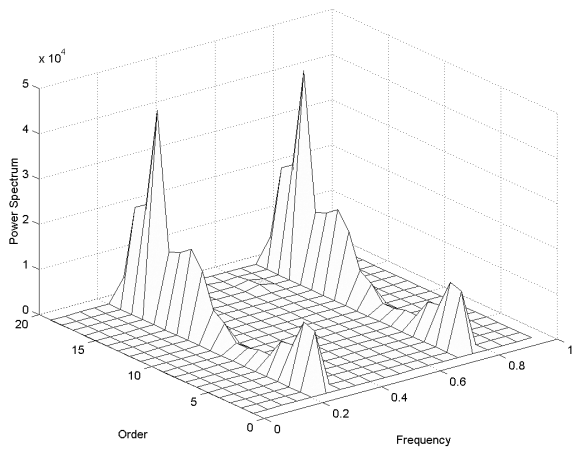

Figure 3

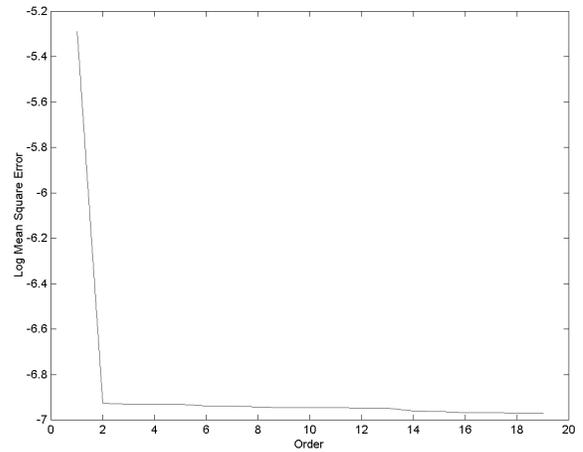

Figure 6

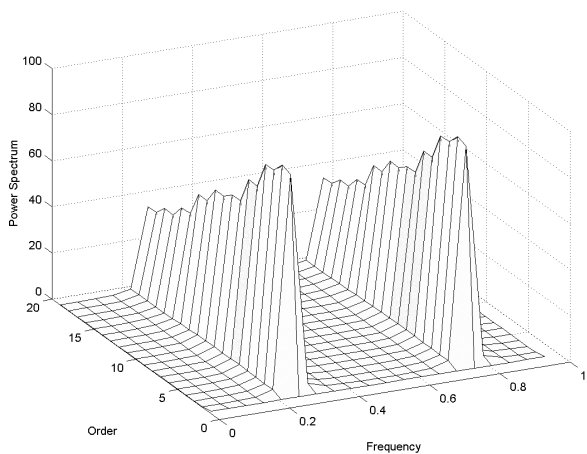

Figure 4

## VI. CONCLUSION

The presented work, establish a theoretical and practical connection between Levinson Algorithm and Burg Algorithm, in 1D and 2D cases. The fact that Unmodified Burg Algorithm fail to minimize the residual error, for certain Order, is of a great importance to multiple application, especially, spectral estimation, where the precision and solution stability is required, and signal compression where error minimization is a priority. In future works we will study the numerical efficiency aspect of the Burg approach relative to the Levinson approach, and present novel ways of applying Levinson and Burg algorithms that do not suffer from Line Splitting phenomena [17].

**R. KANHOUCHE** (M'76–SM'81–F'87) was born in France, Mozelle, on May 8, 1974. He received Engineering Degree in electronics from Damascus University, Damascus Syria, in 1999, and DEA degree (First year Masters) MVA (Mathematiques, Vision, and Automated learning) from Ecole Normale Superieure de Cachan, Cachan, France, in 2001.